\documentclass[12pt,a4paper]{amsart}

\usepackage{amssymb}

\usepackage{amscd}
\usepackage[colorlinks=true, linkcolor=blue]{hyperref}

\def\Bbb{\mathbb}
\def\Cal{\mathcal}

\let\phi\varphi

\newcommand{\x}{\times}

\newcommand{\cD}{\mathcal{D}}
\newcommand{\cN}{\mathcal{N}}

\newcommand{\al}{\alpha}
\newcommand{\be}{\beta}

\newcommand{\de}{\delta}
\newcommand{\ep}{\epsilon}

\newcommand{\la}{\lambda}

\newcommand{\si}{\sigma}

\newcommand{\Ph}{\Phi}

\newcommand{\rpl}                         
{\mbox{$
\begin{picture}(12.7,8)(-.5,-1)
\put(0,0.2){$+$}
\put(4.4,3.1){\oval(8,8)[r]}
\end{picture}$}}

\newcounter{theorem}
\newtheorem{thm}[theorem]{Theorem}
\newtheorem*{thm*}{Theorem \thesubsection}
\newtheorem{lemma}[theorem]{Lemma}

\newtheorem{cor}[theorem]{Corollary}
\newtheorem*{lemma*}{Lemma \thesubsection}
\newtheorem*{prop*}{Proposition \thesubsection}
\newtheorem*{cor*}{Corollary \thesubsection}

\theoremstyle{definition}
\newtheorem{definition}[theorem]{Definition}
\newtheorem*{definition*}{Definition \thesubsection}
\newtheorem{example}[theorem]{Example}
\newtheorem*{example*}{Example \thesubsection}
\theoremstyle{remark}
\newtheorem{remark}[theorem]{Remark}
\newtheorem*{remark*}{Remark \thesubsection}

\def\sideremark#1{\ifvmode\leavevmode\fi\vadjust{\vbox to0pt{\vss
 \hbox to 0pt{\hskip\hsize\hskip1em
 \vbox{\hsize1in\tiny\raggedright\pretolerance10000
  \noindent #1\hfill}\hss}\vbox to8pt{\vfil}\vss}}}%
                        
\begin{document}

\title{Non-holonomic equations for the normal extremals in 
geometric control theory}

\author{A.\ Rod Gover and Jan Slov\'ak}

\address{A.R.G.:Department of Mathematics\\
  The University of Auckland\\
  Private Bag 92019\\
  Auckland 1142\\
  New Zealand\\
J.S.: Department of Mathematics and Statistics\\ Masaryk University, Faculty
of Science\\ Kotlarska 2 \\ 611 37 Brno\\ Czech Republic} 
\email{r.gover@auckland.ac.nz}
\email{slovak@muni.cz}

\begin{abstract} We provide a new and
simple system of equations for the normal sub-Riemannian
geodesics. These use a partial connection that we show is canonically
available, given a choice of complement to the distribution. We also
describe  conditions which, if satisfied, mean that even this
choice of complement is determined canonically, and that this determines
a distinguished connection on the tangent bundle. 

Our approach applies to sub-Riemannian geometry the point of view of
non-holonomic mechanics. The geodesic equations obtained  split into 
mutually driving
horizontal and complementary parts, and the method allows for particular
choices of nice coframes. We illustrate this feature on examples of contact
models with non-constant symbols.

\end{abstract}



\thanks{A.R.G. gratefully acknowledges support from the Royal
Society of New Zealand via Marsden Grant 16-UOA- 051; 
J.S. gratefully acknowledges support from the Grant Agency of the Czech
Republic, grant Nr. GA17-01171S, and the hospitality
of the University of Auckland. The authors are also very grateful to Dmitri
Alekseevsky for numerous very helpful discussions on the topic.}

\maketitle

\pagestyle{myheadings} \markboth{A.\ Rod Gover, Jan Slov\'ak}
{Non-holonomic equations for the normal extremals in geometric control theory}

\section{Introduction}

A large class of problems in control theory deal with optimal control
in $\Bbb R^n$ under linear control and quadratic costs.  A more
geometric formulation leads to the problem of seeking shortest curves
with respect to a Riemannian metric on a finite dimensional manifold,
subject to
certain linear constraints on their velocities.  In the case of no
constraints, we arrive at the Riemannian geodesic curves.  On the
other hand, if the constraints are linear and holonomic, the entire
space foliates to leaves and the controls do not allow the
trajectories to leave these.  The most interesting cases in practice
involve  non-holonomic linear constraints; that is we seek  extremals tangent
to a given distribution $\Cal D$ in the tangent space to the configuration
space. This invokes the
Carnot-Caratheodory metric, and the so-called sub-Riemannian geodesics
are the optimal control curves. In this paper, we 
develop a new geometric approach to treat local problems of the latter
type.

All manifolds will be smooth, connected and finite dimensional, and
all mappings and tensors will also be assumed to be smooth. In
addition to standard smooth affine connections $\nabla$ on manifolds
$M$, we shall also deal with partial connections (and denote these
also by $\nabla$) that directly provide parallel transport only in the
directions of a given distribution $\cD\subset TM$. As a point of notation, typically we
shall use the same notation for bundles and the spaces of their smooth
sections.

As is well known, there are three equivalent approaches to the geodesics on a
Riemannian manifold $(M,g)$: they are locally the shortest curves joining their
points (the variational approach leading to the Euler-Lagrange equations); 
they are the projections of the solutions to the Hamiltonian equations on
the cotangent space $T^*M$ (the Hamiltonian approach minimizing the energy
associated with the curves, i.e.\ the geodesics are the projections of the
flows of the Hamiltonian vector field corresponding to the quadratic 
Hamiltonian $H(p)= \frac12g^{-1}(p,p)$); finally,  they are also the 
curves with autoparallel tangents with respect to the Levi-Civita
connection on $M$.
 
We pass now to the sub-Riemannian situation. In this case we have a metric $g$
defined only on a linear sub-bundle $\cD\subset M$, we
consider the curves that are everywhere tangent to $\cD$, the so
called horizontal curves, and we seek the length minimizers among
them.  The Carnot-Caratheodory (or sub-Riemannian) distance $d(x,y)$
is defined as the infimum of the length of horizontal curves and the
celebrated Chow-Rashevskii theorem says, that this is indeed a metric,
provided the distribution $\cD$ is
bracket generating, i.e. the so called H\"ormander's condition holds
true (cf. \cite{ABBbook, Mont}).
Moreover this metric is topologically equivalent to the Riemannian metric 
for any extension of $g$ to the entire space.  As explained carefully in
\cite{VF10}, the three approaches to distinguished curves mentioned
above are all very interesting from the point of view of non-holonomic
mechanics and more widely, but they provide completely different
concepts in the sub-Riemannian context.

It would seem that the horizontal autoparallel curves have been the most
important ones from the point of view of non-holonomic mechanics. 
This has led to their description in terms of partial connections, 
which encode the geometry in
question with the help of a chosen complement $\cD^\perp$ to $\cD$ in
the tangent bundle. This idea has been known for many decades, and perhaps
goes back to Schouten. See \cite{EK} for a quite detailed account on the
history, and this includes reference to \cite{Sch28}. Numerous authors have discussed the
conditions under which some of the horizontal parallel curves will
happen to be simultaneously the length minimizers, see
e.g. \cite{BBRR, DiVe, VF10}, and also  various notions of curvature of
sub-Riemannian geometry have been discussed, see e.g. \cite{ABBbook,
  AlMeSl, Hla}.

The goal of this paper is different. We derive a new and practical
differential system for the geodesics by exploiting the intuitive
description of the sub-Riemannian minimizers, on a non-holonomic
Riemannian structure $(M,g_\de,\cD,\cD^\perp)$, that arises by
considering the impact of rescaling the costs of the complementary
components of the velocities toward infinity (see the introductory
explanations in the book \cite{Mont}).  Technically and more precisely
this goes as follows. First, we fix an extension of the metric on $\cD$
to the entire tangent bundle, then rescale the metric on $\cD^\perp$
by constants $\ep=\frac 1\de$ and consider $\de\to 0$. We observe that
there is a particularly nice metric connection $\nabla$ respecting the
splitting $TM=\cD\oplus\cD^\perp$ which does not depend on the
parameter $\de$. See Theorem \ref{key3}. Next, we write the
variational equations for the Riemannian geodesics (with $\de>0$) in
terms of the distinguished connection $\nabla$ and renormalize the
complementary components of velocity so that they do not vanish as
$\de\to 0$. In the limit $\de\to 0$ we obtain the equations equivalent
to the Hamiltonian equations for the sub-Riemannian minimizers, thus
providing all normal extremals of the original geometric control
theory problem.

The solutions live in nice geometric bundles that are reminiscent of
the classical tractor bundles and connections originating in the work
of Tracy Thomas nearly one hundred years ago, cf. \cite{BEG}.  We
explain these links in the next section. The technical core of the
paper is explained in section 3 and it is here that we develop the
novel approach to equations for the normal extremals of the control
problems. The coefficients in these equations are given by several
torsion components of the above mentioned connection $\nabla$ and
there is an explicit analytic expression for them, based on the choice
of frames for $\cD$ and $\cD^\perp$. This allows us to discuss the
right choice of the complement $\cD^\perp$ in several examples in
section 4, and to provide some further results giving canonical
choices in section 5.

%
%

In summary our main results are as follows. Theorem \ref{key3} provides a new and
simple system of equations for the normal sub-Riemannian geodesics. In
these equations $u^i$ is a section of the distribution $\cD$ and
$\nu_a$ is a section of the distribution annihilator in $T^*M$. (The
quantities $u^i$ and $\nu_a$ may be viewed as components in suitably adapted
non-holonomic frames.) The
tensor field $L^a{}_{ij}$ is the distribution Levi-form. Then $\nabla$
is a partial connection on the system $(u^i,~\nu_a)$ that has linked
to it some torsion quantities that are the tensors $T^{j}{}_{ak}$ and
$T^b{}_{ak}$ in the system. These latter objects are obtained in the
next main result, which is Theorem \ref{thm_connection}. That Theorem
shows that a choice of complement, in $TM$, to the distribution $\cD$
determines a canonical connection on $TM$ (and hence also determines
the partial connection $\nabla$ and the torsions mentioned). Finally
in Theorem \ref{key1} we show that in the case that the Levi-form is
surjective as a map from $\cD\times \cD\to TM/\cD$ and is also injective in
the sense of equation (\ref{helps}) then in fact no choices are
required: the connection of Theorem \ref{thm_connection} is then
canonical, and hence so also are all the mentioned quantities in the
equations of Theorem \ref{key3}. Theorem \ref{key1} then gives some
general settings where the required surjectivity and injectivity hold. In
particular it is seen there that constancy of the sub-Riemannian symbol is not
necessary in order for the premise of Theorem \ref{key1} to hold.
Canonical metric extensions and choices of complements are discussed in earlier
papers in various contexts, see e.g. \cite{DiVe, Hla}. Whenever a
canonical complement is available, this can be used to remove the choice in 
Theorem \ref{thm_connection}  and hence the connection in Theorem
\ref{thm_connection} and system in Theorem \ref{key3} are completely
canonical.

\section{Sub-Riemannian and non-holonomic Riemannian structures}

A sub-Riemannian geometry consists of a smooth manifold $M$ equipped
with a completely non-integrable distribution 
$\cD\subset TM$, and on $\cD$ a metric, 
i.e. a symmetric positive definite 
bilinear form. Locally, in the optimal control theory problems 
it can be given
by a positive definite matrix valued function $g_{ij}(x)$ on $\Bbb R^n$
together with the suitable
constraint expressed by a field of linear functions
$F(x):\Bbb R^n\to \Bbb R^k$ whose kernel is $\cD\subset \Bbb R^n\x \Bbb
R^n$.   
However the minimal data, that involves no additional choices, can be
described succinctly as follows.
\begin{definition}
A sub-Riemannian geometry is a smooth $n$-manifold $M$ equipped with a
symmetric bilinear (2,0)-tensor $h$ with the properties  that:
(i) the image
$\cD$ of
\begin{equation}\label{key}
h:T^*M \to TM 
\end{equation}
is a  distribution of constant rank $k$ that is bracket generating; and (ii)
the metric on $\cD$ induced by $h$ is positive definite.
\end{definition}

The metric on $\cD$ (determined by the contravariant tensor $h$)  
will also be denoted $h$; by dint of context this should not cause confusion. 

\newcommand{\cK}{\mathcal{K}}
\newcommand{\ccD}{\mathcal{Q}}
Thus on a sub-Riemannian manifold we have, from (\ref{key}), the canonical short exact sequences
\begin{equation}\label{f1}
0\to \cK\to T^*M\stackrel{h}{\to} \cD\to 0,
\end{equation}
and 
\begin{equation}\label{f2}
0\to \cD \to TM \stackrel{q}{\to} \ccD \to 0.
\end{equation}
There is also the $\ccD$-valued Levi-form defined by projecting the Lie bracket of
vector fields in $\cD$
\begin{equation}\label{Levi}
L: \cD\times \cD \to \ccD.
\end{equation}

Note that for any splitting $s:\ccD\to TM$ of (\ref{f2}) we have
$q\circ s=id_\ccD$ and so we can write an element $v\in TM$ as a vector 
direct  sum $v= (s\circ q) (v) + (v-(s\circ q) (v)) = \si + u$ that we represent by a pair
$$
\left( \begin{array}{c} \si^a \\
u^i\end{array}\right), \qquad \mbox{where} \qquad  \si\in \ccD, \, u\in \cD.
$$
Each such  splitting $s$ of   (\ref{f2}) is equivalent to a splitting of  (\ref{f1}) and so we can similarly write an element of
$T^*M$ as a pair
$$
\left( \begin{array}{c} u^i \\
\nu_a \end{array}\right),  \qquad \mbox{where} \qquad u^i\in \cD, \, \nu_a \in \cK,
$$ 
or sometimes $(u^i, \, \nu_a) $. Via the
inclusion $\cK\subset T^*M $ of (\ref{f1}) we can view $\nu$ as an element of $TM$
that satisfies $h(\nu, \ )=0$. 

A change of splitting from $s$ to another $\widehat{s}:\ccD\to TM$
satisfies $q\circ (\widehat{s}-s)=0$ and so $(\widehat{s}-s)$ may be
naturally identified with a bundle map $ f: \ccD \to \cD.  $ 
 Thus under such
a change of splitting the pair representing an element of $TM$  transforms according to 
$$
[TM]_s\ni[v]_s=\left( \begin{array}{c} \si^a \\
u^i\end{array}\right)_s \mapsto \left( \begin{array}{c} \widehat{\si}^a \\
\widehat{u}^i\end{array}\right)_{\widehat{s}}= \left( \begin{array}{c} \si^a \\
u^i - f^i_a\si^a\end{array}\right)_{\widehat{s}}=[v]_{\widehat{s}}\in
[TM]_{\widehat{s}}.
$$ 
In the entire paper, we use the convention of summation over repeated upper
and lower indices.

Similarly for an element of $T^*M$, or sections representing a
section of $T^*M$, we have
\begin{equation}\label{T*change}
\left( \begin{array}{c} u^i \\
\nu_a \end{array}\right) \mapsto \left( \begin{array}{c} u^i \\
\nu_a + f^i_a u_i\end{array}\right) \qquad \mbox{where}\qquad u_i=h_{ij}u^j, 
\end{equation}
and we have simplified the notation. 

If we assume that $\cD$ is bracket generating, then a choice of a
splitting $s$ is equivalent to defining a non-holonomic Riemannian
geometry with the decomposition $TM=\cD\oplus \cD^\perp$, where
$\cD^\perp=s(\ccD)$. In the reverse direction the implication is
obvious.  For the other we proceed as follows.  First note that the
metric $h$ on $\cD$ canonically induces a metric on all tensor bundles
built from $\cD$. Then, second, that (since the distribution is
completely non-integrable) the entire $TM$ may be recovered as a
surjective projection from a suitable such tensor bundle. Thus $TM$ can be identified with the orthogonal complement to
the kernel of the latter projection and it is equipped with the
canonical extension $g$ of the metric $h$ to $TM =
\cD\oplus\cD^\perp$.

\begin{remark} 
Once we are given any non-holonomic Riemannian manifold
$(M,g,\cD,\cD^\perp)$, we can consider a smooth family of the generalized 
sub-Riemannian structures in the sense of \cite{ABBbook}. The general
concept is based on a Riemannian vector bundle $E\to M$ equipped with the
linear control function $\Ph:E\to TM$. Every Lipschitz curve $c(t)$ in $M$ 
allows for the
unique optimal control covering the velocity curve $\dot c(t)$ (defined almost
everywhere) via $\Ph$, i.e. for each $\dot c(t)$ we choose the preimage in
$\Ph^{-1}(\dot c(t))$ of minimal length. 
Thus we can measure the length of all Lipschitz curves and 
the notion of length minimizing curves is a well-defined concept.

In our case of a non-holonomic Riemannian structure $(M,g,\cD,\cD^\perp)$, 
we choose $E=TM$ and for a non-negative real parameter $\al$ we define 
\begin{equation}\label{generalized-sub-Riemannian}
\Ph_\al = \begin{cases} \operatorname{id}_{\cD} &\mbox{on $\cD$} \\
\al\operatorname{id}_{\cD^\perp} &\mbox{on
$\cD^\perp$.}
\end{cases}
\end{equation}
In particular $\Ph_1$ is the identity on $TM$ and then (for length calculations)
with positive $\al$ 
approaching zero we charge each of the $\cD^\perp$ components of the 
velocities $\dot c(t)$ by a $1/\al$ multiple of its original size
with respect to $g$.
At the  $\al=0$ limit we obtain the original 
sub-Riemannian geometry
and $\Ph$ is the orthogonal projection onto $\cD$. It is well known from the
general theory that the behavior of the generalized sub-Riemannian geometry
is regular 
with respect to such smooth deformations, see
\cite{ABBbook} for details.

If we write $\be=1/\al$ for
the positive values of $\al$, then the generalized sub-Riemannian geometry
with $\Ph_\al$ corresponds to non-holonomic Riemannian geometry with the
original metric $g$ on $\cD$, while the metric on $\cD^\perp$ is modified to
$\be^2 g$.
\end{remark}

\section{Distinguished metric connections and normal minimizers}

Let us fix a non-holonomic Riemannian geometry
$(M,g,\cD,\cD^\perp)$. Similar to the earlier works on non-holonomic
Riemannian geometry, our aim is to express the minimizers of the
corresponding sub-Riemannian geometry by means of a special partial
connection $\bar\nabla$ that facilitates differentiation along the
horizontal curves in $\cD$. Since the sub-Riemannian minimizers are
determined by higher order derivatives, we shall have to couple such
curves with some auxiliary parameters in order to reach all of the
minimizers.

We consider the Levi-Civita connection $D$ of 
the metric $g$ on $M$ and define 
\begin{equation}\label{Schouten_D}
\nabla_X Y = (D_X Y)_{\cD}
,\end{equation}
for all vector fields $Y\in \cD$, $X\in TM$, and where 
the subscript denotes the orthogonal projection to $\cD$. 
Symmetrically, we extend the definition to
\begin{equation}\label{Schouten_D^perp}
\nabla_X Y = (D_X Y)_{\cD^\perp}
\end{equation} 
for all vector fields $Y\in \cD^\perp$, $X\in TM$. 

Clearly, the two formulae (\ref{Schouten_D}) and
(\ref{Schouten_D^perp}) together define a metric connection $\nabla$
on $M$ preserving both $\cD$ and $\cD^\perp$. We shall refer to it as
the \emph{Schouten connection} on $(M,g,\cD,\cD^\perp)$. We shall
write $\bar\nabla$ for the partial connection differentiating only in
the directions $X\in\cD$.
The splitting of the tangent space also defines a
torsion of the partial
connection $\bar\nabla$ restricted to $X, Y\in\cD$,
\begin{equation}\label{connection_D}
T^{\cD}_{\cD\cD}(X,Y) = \bar\nabla_XY-\bar\nabla_YX - [X,Y]_{\cD}
,\end{equation}
which vanishes since it coincides with the orthogonal projection of the (zero) 
torsion of $D$ as a connection on $TM$.

More generally, we shall write 
$T^A_{BC}$ for the components of the torsion of $\nabla$, 
$$
T^A_{BC} :B\x C \to A
,$$
where the letters $A$, $B$, $C$ stand for either $\cD$ or $\cD^\perp$.

\begin{thm}\label{thm_connection} Given a sub-Riemannian geometry $(M,\cD,h)$,
  let $g$ be a Riemannian metric on $TM$ that restricts to $h$ on $\Cal D$ and
write $\Cal D^\perp$ for the orthogonal complement of $\Cal D$. Then there
is the unique metric connection $\nabla$ on $TM$
such that both $\Cal D$ and $\Cal D^\perp$ 
are preserved, and 
\begin{gather}
T^{\Cal D}_{\Cal D\Cal D}= 0\label{partial_torsion}
\\
T^{\Cal D^\perp}_{\Cal D^\perp\Cal D^\perp}= 0\label{partial_torsion_perp}
\\
T^{\Cal D^\perp}_{\Cal D\Cal D^\perp}\mbox{\ is symmetric 
with respect to $g_{|\Cal D^\perp}$} \label{mixed_torsion}
\\
T^{\Cal D}_{\Cal D\Cal D^\perp} \mbox{\ is symmetric with respect to $g_{|\Cal D}$} 
\label{mixed_torsion_perp}.
\end{gather}
\end{thm}

\begin{proof}
  For the purposes of the proof we will denote the claimed new
  connection by $\tilde{\nabla}$, this is to simplify notation in the
  discussion.
     We shall see that the required connection $\tilde{\nabla}$ is a straightforward
modification of the Schouten connection introduced in
\eqref{Schouten_D} and \eqref{Schouten_D^perp}.  By its definition, the
Schouten connection $\nabla$ preserves $g$, $\cD$ and $\cD^\perp$ (i.e.\ these are parallel) 
and two components of its torsion vanish, namely those specified  in
\eqref{partial_torsion}, \eqref{partial_torsion_perp}.

Let us now consider any connection $\nabla$ on $TM$ which leaves the
metric $g$ parallel and write $\langle \ ,\ \rangle$ as a shorthand
for the scalar product $g( \ , \ )$ on $TM$.  We recall the standard defining
equation for the torsion of $\nabla$:
\begin{equation}\label{std_torsion}
T(X,Y) = \nabla_XY-\nabla_YX-[X,Y] .
\end{equation}
For arbitrary vector fields $X,Y,Z$ 
this leads to
\begin{equation}\label{full_nabla_eq}
\begin{aligned}
X\langle Y,Z\rangle &+ Y\langle Z,X\rangle - Z\langle X,Y\rangle = \langle
\nabla_XY,Z\rangle +\langle Y,\nabla_XZ\rangle 
\\ &\qquad\qquad + \langle
\nabla_YZ,X\rangle +\langle Z,\nabla_YX\rangle -\langle
\nabla_ZX,Y\rangle -\langle X,\nabla_ZY\rangle
\\
&= 2\langle \nabla_XY,Z\rangle - \langle Z,[X,Y]+T(X,Y)\rangle - \langle
Y,[Z,X]+T(Z,X)\rangle 
\\
&\qquad\qquad + \langle X,[Y,Z]+T(Y,Z)\rangle
,\end{aligned}
\end{equation}
since $\nabla $ is metric.
We are going to restrict the latter two equations, i.e.\ (\ref{std_torsion}) and  (\ref{full_nabla_eq}), to special choices of the
arguments to derive necessary conclusions for $\nabla$.

First we  use the preservation of $\cD^\perp$ (by $\nabla$) and
\eqref{std_torsion} to see that
\begin{equation}\label{ppd_torsion}
T^{\cD}_{\cD^\perp\cD^\perp}(X,Y) = -[X,Y]_{\cD}
\end{equation}
for all $X,Y\in \cD^\perp$,
with the projection of the Lie bracket of vector fields to $\cD$ on
the right hand side.

Next, we analyze the equation \eqref{full_nabla_eq} for the same
choice of components for $X,Y,Z$, i.e. $X,Y\in \cD^\perp$ and $Z\in \cD$, 
and use the last result (\ref{ppd_torsion}) to obtain
\begin{align*}
-Z\langle X,Y\rangle &= -\langle Z,[X,Y]\rangle - \langle
Y,[Z,X]\rangle +\langle X,[Y,Z]\rangle \\&\qquad - \langle
Z,-[X,Y]\rangle -\langle Y,T(Z,X)\rangle +\langle X, T(Y,Z)\rangle
.\end{align*} 
From this we see that the symmetric part of 
$T^{\cD^\perp}_{\cD\cD^\perp}$  is completely
determined:
\begin{equation}\label{dpd_torsion}
\langle X,T(Y,Z)\rangle + \langle Y,T(X,Z)\rangle= -Z\langle X,Y\rangle + \langle Y,[Z,X]\rangle +
\langle X,[Z,Y]\rangle
.\end{equation}
In particular, the choice $X=Y$ gives
\begin{equation}\label{diag_dpd}
\langle X, T(Z,X)\rangle = \frac12 Z\|X\|^2 +\langle X,[X,Z]\rangle 
.\end{equation}

Let us check what is the $T^{\cD^\perp}_{\cD\cD^\perp}$ component of
the torsion of the Schouten connection $\nabla$.
\begin{align*}
T^{\cD^\perp}_{\cD\cD^\perp}(Z,X) &= \nabla_ZX-[Z,X]_{\cD^\perp}
\\
&=(D_ZX-[Z,X]-D_XZ+D_XZ)_{\cD^\perp} = (D_XZ)_{D^\perp}
\end{align*}
since the torsion of $D$ vanishes. Its symmetric part is given by
$$
\langle X,T(Z,X)\rangle = \langle X,D_XZ\rangle = \langle X,D_ZX+[X,Z]\rangle
,$$ 
as predicted by \eqref{diag_dpd}, while 
the antisymmetric part is
\begin{equation}\label{schouten_torsion}
\begin{aligned}
\langle Y,T^{\cD^\perp}_{\cD\cD^\perp}(Z,X)\rangle
&- \langle X,T^{\cD^\perp}_{\cD\cD^\perp}(Z,Y) \rangle
= \langle Y,D_XZ\rangle - \langle X,D_YZ\rangle 
\\
&= \langle D_YX,Z\rangle -
\langle D_XY,Z\rangle
=\langle [Y,X],Z\rangle
.\end{aligned}
\end{equation}

Thus, in order to satisfy \eqref{mixed_torsion} and
\eqref{mixed_torsion_perp}, we have to deform $\nabla$. The only option is to
posit that $ \tilde T_{\cD\cD^\perp}^{\cD^\perp}$, defined by
\begin{equation}\label{pd_derivative}
\tilde\nabla_XY = \tilde T_{\cD\cD^\perp}^{\cD^\perp}(X,Y) +
[X,Y]_{\cD^\perp}
\end{equation}
with $X\in\cD$ and $Y\in\cD^\perp$, is symmetric and that this
symmetric torsion component is defined by \eqref{diag_dpd}.  We keep
the derivative $\nabla_XY$ unchanged for $X$, $Y\in\cD$.
If we exchange
the roles of $\cD$ and $\cD^\perp$ in the above considerations, we
obtain the relevant formulae for the derivatives $\tilde\nabla_XY$ of
the deformed Schouten connection $\tilde \nabla$.  It remains to check
that the connection $\tilde \nabla$ satisfies all the conditions of
the theorem.

The four components $\tilde T_{\cD\cD}^{\cD^\perp}$, 
$\tilde T_{\cD^\perp\cD^\perp}^{\cD}$, $\tilde T_{\cD\cD}^{\cD}$, and 
$\tilde T_{\cD^\perp\cD^\perp}^{\cD^\perp}$ have not changed. The remaining 
two components are symmetric, by definition. 

Finally, we 
have to check that the deformed connection $\tilde \nabla$ preserves 
the metric $g$, i.e. we need to check 
$Z\langle X,Y \rangle=\tilde\nabla_Z\langle X,Y \rangle = \langle
\tilde\nabla_ZX,Y\rangle + \langle X,\tilde\nabla_ZY\rangle$. Clearly this
has to be checked for $X$, $Y$ from the same component, say $\cD$. If $Z$ is
from the same one, nothing is changed compared to the Schouten connection 
and so the condition holds true. If
$Z\in\cD^\perp$, then exploiting the fact that
$\nabla$ preserves
both $\cD$ and $\cD^\perp$ and using (\ref{pd_derivative}), we may rewrite this condition as
$$
\tilde\nabla_Z\langle X, Y\rangle = \langle \tilde T(Z,X),Y\rangle + \langle
\tilde T(Z,Y),X\rangle + \langle [Z,X],Y\rangle + \langle [Z,Y],X\rangle
.$$ 
By the definition of $\tilde\nabla$, its mixed 
torsion components differ from the Schouten connection only in their
antisymmetric parts and, thus,
$$
\tilde\nabla_Z\langle X, Y\rangle = \langle T(Z,X),Y\rangle + \langle
T(Z,Y),X\rangle +\langle [Z,X],Y\rangle + \langle [Z,Y],X\rangle = \nabla_Z\langle X,Y\rangle
,$$
as requested. 

It is well known, and easily verified, the there is exactly one metric
connection on $M$ for each prescription of its torsion. Thus, our
connection $\tilde \nabla$ is the unique one satisfying the
assumptions of the theorem.
\end{proof}

For later use, let us notice that the derivatives $\nabla_XY$ for both $X$,
$Y$ from the same component are given by the usual formula known for the
Levi-Civita connections (with $Z$ from the same component):
\begin{equation}\label{levi_civita_eq}
\begin{aligned}
2\langle \nabla_XY,Z\rangle &= 
X\langle Y,Z\rangle + Y\langle Z,X\rangle - Z\langle X,Y\rangle
\\
&\qquad  
+ \langle Z,[X,Y]\rangle + \langle
Y,[Z,X]\rangle 
- \langle X,[Y,Z]\rangle
.\end{aligned}
\end{equation}

The metric connection $\nabla$ of Theorem \ref{thm_connection}, as
constructed explicitly in the proof, has a nice and useful
property under constant rescalings of the metrics involved:

\begin{cor}\label{three}
Let $\tilde g$ be another metric on $\Cal D\oplus \Cal D^\perp$ that
 differs from $g$ only by constant rescaling on each of the orthogonal complements
$\Cal D$ and $\Cal D^\perp$. Then the corresponding unique connection
$\tilde \nabla$ from  Theorem \ref{thm_connection} remains
unchanged. 
\end{cor}

\begin{proof}
We have to analyze the relevant formulae in our proof of Theorem
\ref{thm_connection}. The equation \eqref{levi_civita_eq} involves scalar
products of one of the orthogonal complements only, thus the definition of
this part of the connection is independent on the constant rescaling (as
well known from the Riemannian geometry). 

Similarly, the derivative $\tilde \nabla_XY$ for fields in different
components is given in \eqref{pd_derivative} where the metric enters via 
\eqref{diag_dpd}. Again, only metric on one of the components enters and the
formulae do not see any constant rescaling. 
\end{proof}

\begin{remark}\label{summary}
It is useful to summarize the argument leading to the connection
above. The metric $g$ on $TM$ provides the orthogonal complement to
$\cD$ in $TM$, thus we have $TM=\cD\oplus \cD^\perp$ and we request these
subbundles to be parallel. Then the torsions
$T^{\cD^\perp\cD^\perp\cD}$ and $T^{\cD\cD \cD^\perp}$ are
determined (see (\ref{ppd_torsion}) and its $\cD \leftrightarrow
\cD^\perp$ swap), as are the symmetric parts of
$T^{\cD^\perp\cD\cD^\perp}$ and $T^{\cD\cD^\perp \cD}$, via
(\ref{dpd_torsion}) and its $\cD \leftrightarrow \cD^\perp$ swap. The
connection is then uniquely determined by setting all remaining
torsion components to zero.

Finally in this remark, let us note here that the Schouten connection does not enjoy the useful
Corollary \ref{three} property that $\tilde{\nabla}$ has, see \eqref{schouten_torsion}.
\end{remark}

Our next goal is to find equations for the horizontal metric geodesics
(metric minimizers among the curves tangent to $\Cal D$) by means of
our distinguished connection $\nabla$. (Note we are now dropping the
temporarily introduced tilde.)

As before, 
let us fix some extension metric $g$ of the given sub-Riemannian metric $h$,
write $TM=\Cal D\oplus \Cal D^\perp$, and consider the family of metrics
$g^\ep$ with $g^\ep_{|\Cal D}=g$ and $g^\ep_{|\Cal D^\perp}=\ep g$. 
Notice that the corresponding linear connection $\nabla^\ep$ remains
unchanged, thus we shall use the same symbol $\nabla$ for all of them.

At the same time, the Riemannian geodesics $c^\ep(t)$ joining the same
points $x_0, x_1\in M$ will depend on $\ep$
heavily. With growing $\ep\to\infty$, the $\Cal D^\perp$ directions
on the geodesics are charged $\sqrt\ep$ times more and thus they become
horizontal curves in the limit (if such a limit exists).

We want to understand the geodesic equation for the metric minimizers of
$g^\ep$ in term of $\nabla$ and its torsion. Let us write $D^\ep$ for the
Levi Civita connection of $g^\ep$ and let $A^\ep:TM\otimes TM \to TM$ be the contorsion tensor
defined by
\begin{equation*}
D^\ep_XY = \nabla_XY + A^\ep(X,Y)
.\end{equation*}
It is well known that given a Riemannian metric and choosing any torsion
tensor, there will be exactly one metric connection with the chosen torsion.
Thus the contorsion tensor $A^\ep$ is uniquely determined by the
torsion $T$ of our connection $\nabla$. Moreover, the $A^\ep$ must be
antisymmetric with respect to the metric $g^\ep$ since both $D^\ep$ and
$\nabla$ preserve the metric.

We shall work in local non-holonomic frames spanning $\Cal D$ and $\Cal
D^\perp$ and we shall use abstract indices $i,j,k,\dots$ and $a,b,c,\dots$
in relation to $\Cal D$ and $\Cal D^\perp$,
respectively. In particular, let us write $u=u^i+u^a$ for the tangent curve
$u=\dot c$, and $\nabla=\nabla_i+\nabla_a$ for the connection. 
Similarly, our fixed metrics are the products of $g_{ij}$ and
$\ep g_{ab}$, while using the decomposition one has that the torsion is the sum of components
$$
T^i{}_{jk}+T^i{}_{ja} + T^i{}_{ab} + 
T^a{}_{jk}+T^a{}_{jb} + T^a{}_{bc} .
$$ In fact the first and the last components vanish for our
connection, cf.\ Theorem \ref{thm_connection}.

\begin{lemma}
The variational equations $D^\ep_u u=0$ 
for the tangent curves $u=\dot c{}^\ep$ of the $g^\ep$ critical curves $c^\ep$ are
\begin{equation}\label{var_via_nabla}
\begin{aligned}
0 &= g_{ij}u^k\nabla_k u^j + g_{ij}u^a\nabla_a u^j + g_{kj}u^k T^j{}_{ia}u^a 
+ \ep g_{ab}u^a T^b{}_{ic}u^c + \ep g_{ab}u^a T^b{}_{ik}u^k
\\
0 &= \ep g_{ab} u^k\nabla_ku^b + \ep g_{ab} u^c\nabla_c u^b 
+ g_{ij}u^i T^j{}_{ab}u^b + g_{ij}u^i T^j{}_{ak} u^k  
 + \ep g_{cb}u^b T^c{}_{ak}u^k.
\end{aligned} 
\end{equation}
\end{lemma}

\begin{proof}
Let us first recall some well-known facts about the critical
curves with respect to the Riemannian length functional.
For each Riemannian metric $g$, the critical points of the variations 
of a curve $c(t)$ with fixed points $c(0)$ and $c(1)$, 
parametrized by a constant multiple of
length, are given by the equation $0=g(\be,D_uu)$, where $\be$ is the
derivative of the variation. 
Thus we are interested in the equations $0= g^\ep(\be,\nabla_u u + A^\ep(u,u))$
for arbitrary values of $\be$,
but we need them written down 
explicitly in terms of the torsion $T$ of $\nabla$.

The defining equation of torsion says (recalling the torsion of $D^\ep=0$)
$$
 T(X,Y) = A^\ep(Y,X) - A^\ep(X,Y) 
.$$
Thus, $\langle Z, T(X,Y)\rangle^\ep = \langle Z, A(Y,X) - A(X,Y)\rangle^\ep$.
If we subtract the same expressions with cyclic permutations of $X,Y,Z$, we
arrive at
\begin{equation*}
2\langle X,A(Y,Z)\rangle^\ep = \langle Z,T(X,Y)\rangle^\ep - 
\langle X,T(Y,Z)\rangle^\ep - \langle Y,T(Z,X)\rangle^\ep,
\end{equation*}
where we exploited the antisymmetry of the contorsion tensor $\langle
\cdot ,A( \cdot,Z)\rangle $.  
The expression
we are interested in is
\begin{equation*}
g^\ep(\be,\nabla_uu+A(u,u)) = g^\ep(\be,\nabla_uu) + g^\ep(u,T(\be,u))
.\end{equation*}

Finally, we expand the  expression on the right hand side in terms of the components of
$\be$, $u$ and $T$. We arrive at:
\begin{multline*}
g_{ij}\be^i(u^k\nabla_ku^j +u^a\nabla_au^j)  
+ \ep g_{ab}\be^a(u^k\nabla_ku^b+ u^c\nabla_cu^b)  
\\
+ g_{ij}u^i(T^j{}_{ka}\be^ku^a 
+ T^j{}_{ak}\be^au^k + T^j{}_{ab}\be^au^b) \\
+ \ep g_{ab}u^a(T^b{}_{k\ell}\be^ku^\ell + T^b{}_{kc}\be^ku^c +
T^b{}_{ck}\be^cu^k).
\end{multline*}
Collecting the terms with
$\be^i$ and $\be^a$ separately, we establish the independently 
vanishing sets of
equations, exactly as in the proposition of our lemma. 
\end{proof}

Now we are in position to analyze the limit behavior of the metric
minimizers. In order to understand the equations better, we shall rename the
$\Cal D^\perp$ component $u^a$ as 
$$
u^a = \frac1\ep \nu^a
.$$
Under this change, writing $\de=1/\ep$, 
the equations \eqref{var_via_nabla} become
\begin{equation}\label{eq_with_nu}
\begin{aligned}
0 &= g_{ij}u^k\nabla_k u^j + \de g_{ij}\nu^a\nabla_a u^j + 
\de g_{kj}u^k T^j{}_{ia}\nu^a 
+ \de g_{ab}\nu^a T^b{}_{ic}\nu^c + g_{ab}\nu^a T^b{}_{ik}u^k
\\
0 &= g_{ab} u^k\nabla_k\nu^b + \de g_{ab} \nu^c\nabla_c \nu^b  
+ \de g_{ij}u^i T^j{}_{ab}\nu^b + g_{ij}u^i T^j{}_{ak} u^k  
 + g_{cb}\nu^b T^c{}_{ak}u^k.
\end{aligned} 
\end{equation}
This is a smoothly parametrized system of differential equations and we are
most interested in the limit for $\de=0$. This is the limit case of the
deformed sub-Riemannian geometry in \eqref{generalized-sub-Riemannian}, and
all other postitive values of $\de$ describe the geodesics of regular Riemannian
metrics. 

Using our original fixed metric
$g$ to lower indices, we may rewrite the limit equations with $\de=0$
\begin{equation}\label{limit_eq}
\begin{aligned}
0 &= g_{ij}u^k\nabla_k u^j + g_{ab}\nu^a T^b{}_{ik}u^k
\\
0 &= g_{ab} u^k\nabla_k\nu^b + g_{ij}u^i T^j{}_{ak} u^k  
+ g_{cb}\nu^b T^c{}_{ak}u^k
\end{aligned} 
\end{equation}
as equations coupling the components $(u^i)\in \Cal D$ with $(\nu_a)$
in the annihilator of $\Cal D$ in $T^*M$ which we shall again denote
as $\Cal D^\perp$ (and identify with $\cD^\perp\subset TM$ via the
metric $g_{ab}$).  Thus everything gets an intrinsic meaning from the
point of view of the sub-Riemannian geometry, except the torsions
which reflect our choice of the complement to $\Cal D$ and the metric
on it, i.e.  the non-holonomic Riemannian extension of
$(M,\cD,h)$. Moreover, notice that \eqref{ppd_torsion} (with $\Cal D$ and $\Cal
D^\perp$ swapped) reveals that
$T^a{}_{ik}$ actually coincides with the Levi form $L^a{}_{ik}$.  

\begin{lemma}\label{key3-lemma}
Projections $c(t)\in M$ of the solutions $v(t)=(u^i(t),\nu^a(t))$ 
to the equations \eqref{limit_eq} are horizontal curves parametrized by
constant speed.
\end{lemma}
\begin{proof}
Consider a fixed value $(u^i(0),\nu^a(0))\in T_xM$ and 
write $v(t,\de)=(u^i_\de(t),\nu^a_\de(t))$ for the solutions of
equations \eqref{eq_with_nu} with $\de\ge 0$ and the common initial conditions
$(u^i(0),\nu^a(0))$. Since the system of equations is smooth, the
mapping $v(t,\de)$ will be smooth too. In particular, the norm $\|v(t,\de)\|$
with respect to the metric $g$ 
will be 
bounded on compact subsets and therefore the same must be
true for the norm of its $\cD^\perp$ component. 

Now, for all nonzero $\de$, these solutions are
Riemannian geodesics with initial velocity $(u^i(0),\sqrt\de\nu^a(0))$. In
particular they are parametrized with constant velocity in the metric 
$g^{1/\de}$, i.e. in the metric $g$ 
$$ \|v_\de(t)\|^2 = \| u^i_\de(t)\|^2 + \de \|\nu^a_\de(t)\|^2= \|
u^i_\de(0)\|^2 + \de\|\nu^a_\de(0)\|^2 .$$
Thus, the norm of the
$\cD^\perp$ component of the velocity of the geodesics, $\sqrt{\de}
\|\nu^a_\de(t)\|$, must converge to zero.  This implies that the
projection of the resulting curve $v_0(t)$ to the manifold is
horizontal with initial velocity $u^i(0)$.

Finally, we look at the parametrization of the $\cD$ component of a solution
$v(t)$. Its norm $\|u(t)\|$ is easily computed from the first
equations of \eqref{limit_eq}. Indeed, we already know that the projection
curve is horizontal, and thus
$$
\tfrac d{dt}\langle u,u\rangle = 2\langle u^k\nabla_ku,u\rangle
= -2 g_{ab}\nu^aT^b{}_{ik}u^ku^i =0
,$$
since the torsion is antisymmetric in the lower indices. Thus, the norm
$\|u(t)\|$ remains constant. 
\end{proof}

Geometrically, we can interpret the lemma as follows. 
For each initial condition of the horizontal velocity $u(0)$ at the point
$x_0\in M$ (the actual velocity of the expected minimizing curve in the
limit), completed by any choice of $\nu(0)\in \cD^\perp$, 
there is a (locally defined) solution to the
system of equations with $\de=0$. The choice of the initial condition
$\nu(0)$ reflects exactly the expected freedom for
sub-Riemannian geodesics with the given initial velocity $u(0)$ at $x_0$.
In terms of the deformation with $\de>0$, the actual 
$\Cal D^\perp$ components of the velocity vector $u(t)=\dot
c(t)$ of the geodesics become negligible for $\de$ close to zero, but the
constantly rescaled values $\nu$ stay of roughly the same size. 

We are ready to prove the key theorem. Recall there is the Levi form
$L^a{}_{ik}$ and the two symmetric torsion components $T^j{}_{ak}$,
$T^b{}_{ak}$ coming from the choice of the non-holonomic Riemannian
extension, see Theorem \ref{thm_connection}.

\begin{thm} \label{key3}
For each set of initial conditions $x\in M$, $u(0)\in\Cal D\subset T_xM$, and 
$\nu(0)\in\Cal D^\perp\subset T_x^*M$, the component $u(t)$ of the 
unique solution of the equations
\begin{equation}\label{limit_eq_lowered}
\begin{aligned}
0 &= u^k\nabla_k u^i + h^{ij}\nu_a L^a{}_{ik}u^k
\\
0 &= u^k\nabla_k\nu_a + g_{ij}u^i T^j{}_{ak} u^k  
+ \nu_b T^b{}_{ak}u^k
\end{aligned} 
\end{equation}
projects to a locally defined normal extremal $c(t)$ of the sub-Riemannian
geometry with $c(0)=x$ and $\dot c(t)=u(t)$.
\end{thm}

\begin{proof}
The two systems of equations \eqref{limit_eq} and \eqref{limit_eq_lowered}
are clearly equivalent and the $\Cal D$ components of solutions
coincide. 

Let us consider a (locally defined) solution $(u(t),\nu(t))$ of \eqref{limit_eq_lowered}
with the given initial conditions $(u(0),\nu(0))$. As discussed already in
the proof of Lemma \ref{key3-lemma}, the projection $c(t)$ of the 
curve $u(t)\in TM$ is the limit of the geodesics $c_\de(t)$ in the metrics 
$g^{1/\de}$ with initial velocities
$\dot c_\de(0)=(u(0),\sqrt\de\nu(0))$ and all of them are solutions to 
the Hamiltonian
equations for the geodesics. The latter Hamiltonian equations on $T^*M$ 
are again
smoothly dependent on the parameter $\sqrt\de$, and their limit case at $\de=0$ provides
the Hamiltonian equations for the normal extremals of the sub-Riemannian
problem. 
Thus our limit $c(t)$ of the geodesics $c_\de(t)$ must be the 
normal extremal $c(t)$.
\end{proof}

\begin{remark}
A few remarks are due here. First let us notice that our construction
of the distinguished connection $\nabla$, of Theorem
\ref{thm_connection}, and the subsequent computations were not
dependent on the assumption that $\Cal D$ is bracket generating. Only
the local existence of the minimizers would not be guaranteed  if
we remove this assumption. In particular, if both $\Cal D$ and $\Cal
D^\perp$ are involutive, then $(M,g)$ is locally a product of two
Riemannian manifolds, all the torsions disappear and our equations
coincide with the standard equations for geodesics. Expanding one of
the metric components by $\ep$ allows one to find the horizontal minimizers
only within the individual leaves of the foliation.



In the theorem, the initial condition for the parameter
$\nu_a\in T^*M$ in the annihilator of $\Cal D$ are linked to the 
initial acceleration of the minimizer in the direction complementary to
$\Cal D$. As expected, this non-trivial 
vertical acceleration is allowed by the bracket
generating condition on $\Cal D$. The coupled equations on $u^i$ and $\nu_a$
determine the unique evolution of this complementary acceleration.
\end{remark}

\section{Examples}

The equations \eqref{limit_eq_lowered}  
for all normal sub-Riemannian geodesics in the main theorem \ref{key3} 
are related to
non-holonomic frames $\langle X_1,\dots, X_n\rangle = \Cal D$ and $\langle
Z_1,\dots,Z_{\ell}\rangle = \Cal D^\perp$, $n+\ell=m$. 

In this section, we compare them 
to the usual systems of $2m$ 1st-order ODEs in holonomic coordinates
$(x^1,\dots,x^n,z^1,\dots,z^\ell)$ on $\Bbb R^m$, proceeding as
follows. 

We express our solution curve $u(t)\in \Cal D\subset T\Bbb R^m$ 
as $u(t)=\al^iX_i$ which 
allows to express the derivatives $\dot x^i$ of the projection
$c(t)=(x^1(t),\dots,x^n(t),z^1(t),\dots,z^\ell(t))$ of $u(t)$ in terms of
the new quantities $\al^i(t)$. The derivatives of the remaining coordinates
$z^a$ are then given by the non-holonomic constraints, with all $\dot x^i$
substituted by the latter expressions with $\al^j$. In this way, we obtain
the $m$ 1st-order equations which are implicitly hidden as the projection
of the solution $u(t)$ of \eqref{limit_eq_lowered} to $M$. 

The further $n$ 1st-order equations on functions $\al^i$ 
are obtained from the first line of \eqref{limit_eq_lowered}, while the
remaining $\ell$ 1st-order equations on the coupled functions in the
expression $\nu(t) = \nu^a(t)Z_a\in\Cal D^\perp$ come from the second line
in \eqref{limit_eq_lowered}. The torsion coefficients are all easily
expressed by means of the identities \eqref{ppd_torsion}, \eqref{diag_dpd},
and their analogues with $\Cal D$ and $\Cal D^\perp$ swapped. Finally, we
have to express the covariant derivatives of $u$ and $\nu$ in the direction
$u(t)$. The first follows from \eqref{levi_civita_eq} since the derivative
restricted to $\Cal D$ is given by the formula for the Levi-Civita
connection. Thus, leaving out terms which are obviously zero, we arrive at
(suppressing the argument $t$)
\begin{equation}\label{duu}
\begin{aligned}
\langle \nabla_{u}u,X_k\rangle &= \dot\al^k + \al^i\al^j\langle 
\nabla_{X_i}X_j, X_k\rangle
\\
&= \dot\al^k + \frac12\al^i\al^j\bigl(\langle X_i,[X_k,X_j]\rangle +\langle
X_j,[X_k,X_i]\rangle\bigr) 
.\end{aligned}
\end{equation}
Finally, the covariant derivative $\nabla_u\nu$ is given by
\eqref{pd_derivative}, where the torsion term appears just with the opposite
sign than in \eqref{limit_eq_lowered} and so only the projection of the
bracket remains, which splits further into 
\begin{equation}\label{dunu}
\al^i\langle [X_i,\nu^aZ_a],Z_b\rangle = \dot \nu^b + \al^i\langle \nu^a[X_i,Z_a],Z_b\rangle
.\end{equation}

We illustrate this procedure on two examples, 
including one with non-constant symbols.
We are choosing the sub-Riemannian metric so that an orthonormal frame of the
horizontal distribution generates directly a reasonable complement, and our
approach then leads to relatively simple equations. This is quite common in
applications and, technically speaking,
the advantage of our approach consists in minimizing the torsions appearing 
in the equations.

\begin{example}[Free 1-step generating distributions]
The generic $n$-dimensional distributions $\cD$ on manifolds $M$ of dimension
$\frac12n(n+1)$ are called free distributions. Picking any local frame
$X_1,\dots,X_n$ of $\cD$, the Lie brackets $Y_{ij}=[X_i,X_j]$, $i<j$ provide
the frame of a complement $\cD^\perp$. This is one of famous Cartan
geometries.

Let us choose the usual frame of the local homogeneous model in standard 
coordinates $(x^1,\dots,x^n,y^{12},\dots,y^{(n-1)n})$ on $\Bbb
R^{\frac12n(n+1)}$,
\begin{equation}\label{D_frame_free}
\begin{aligned}
X_i &= \frac\partial{\partial x^i} - x^{i+1}\frac\partial{\partial y^{i(i+1)}}
- \dots - x^n\frac\partial{\partial y^{in}},
\ 1\le i\le n
\\
Y_{ij} &= [X_i,X_j] = \frac\partial{\partial y^{ij}}
,\end{aligned}
\end{equation} 
and define the metric $g$ to make this an orthonormal basis.

The dual basis on $T^*M$ is obviously 
\begin{equation}\label{free_condition}
\begin{gathered}
dx^1,\dots, dx^n, 
\\
x^2dx^1+dy^{12},\ x^3dx^1+dy^{13},\dots,\ x^jdx^i+dy^{ij},\dots
\end{gathered}
\end{equation}

Let us write down the equations \eqref{limit_eq_lowered} in our coordinates
explicitly. They consist of two layers. 
First, the above choice of $u(t)=\al^iX_i$ leads directly to 
new names to the derivatives $\dot x^i=\al^i$. Next, the
non-holonomic horizontality condition means the tangent vectors are in the
kernel of the forms in the second line of \eqref{free_condition}
and we arrive at the remaining 
$\frac12n(n-1)$ equations from the first set of $\operatorname{dim}M$ 
equations:
\begin{equation}
\begin{gathered}\label{eq1_free}
\dot x^i = \al^i,\quad 1\le i \le n
\\
\dot y^{ij} = -x^j\al^i,\quad 1\le i<j\le n
.\end{gathered}
\end{equation}
\end{example} 
As discussed above, the next set of 
$\operatorname{dim}M$ equations, i.e. the equations of theorem \ref{key3}, 
are expressed in the non-holonomic frames, cf. \eqref{limit_eq_lowered}. 
We discuss them now.

The coordinate description of the Levi
form is clear from the choice of the frame and provides the first set of
equations below. Further, all brackets of the generators of $\Cal D$ are in
$\Cal D^\perp$ and thus the components of the covariant derivative $\nabla_uu$
in \eqref{duu} are expressed by $\dot \al$ only. 
A direct inspection of the
formula \eqref{diag_dpd} reveals that both mixed torsion components vanish
since we are working with an orthonormal frame (i.e. $\|X_i\|=1$) 
and all the brackets
$[X_i,Y_{jk}]$ vanish identically. Thus, the second set of equations is
trivial and we arrive at 
\begin{equation}
\begin{aligned}\label{eq2_free}
\dot\al^i &= \nu_{1i}\al_1 + \dots \nu_{(i-1)i}\al^{i-1}
-\nu_{i(i+1)}\al^{i+1}-\dots - \nu_{in}\al^n,\quad 1\le i \le n
\\
\dot \nu_{ij} &= 0,\quad 1\le i<j\le n.
\end{aligned}
\end{equation}
It is quite straightforward to solve these equations explicitly, just the
general solution formulae are a bit messy. 
The most trivial initial
condition $\nu_{ij}=0$ implies $\al_i$ are arbitrary constants, $x^i$ are 
affine functions in $t$, while the $y^{ij}$ are then quadratic.

On the other hand, if we choose just one of the $\nu$'s as nontrivial constant $C\ne0$, 
the solutions are similar to the lowest dimensional case $n=2$,
which recovers the most classical three dimensional
Heisenberg group example. In this case we deal with 
(writing $z$ instead of $y_{12}$) 
\begin{align*}
X^1 &= \frac\partial{\partial x^1} - x^2\frac\partial{\partial z}
\\
X^2 &=\frac\partial{\partial x^2}
\end{align*}
and the equations \eqref{limit_eq_lowered},
\begin{align*}
\dot x^1 &= \al^1,\quad
\dot x^2 = \al^2
\\
\dot z &= -x^2\al^1
\\
\dot \al^1 &= - \nu \al^2
,\quad
\dot \al^2 = \nu\al^1
\\
\dot \nu &= 0
\end{align*}
have got the solutions (with fixed constant $\nu=K\ne0$ 
and five free parameters $C_1,\dots,C_5$:
\begin{align*}
x^1(t) &= \frac{C_1}K\cos(Kt+C_2)+C_3
\\
x^2(t) &= \frac{C_1}K\sin(Kt+C_2)+C_4
\\
z(t) &=
\frac{C_1{}^2}{2K} t - \frac{C_1C_4}K\cos(Kt+C_2)-\frac{C_1{}^2}{4K}\sin(2Kt+2C_2)
+\frac{C_1^2C_2}{2K}+C_5.
\end{align*}
Although the coordinates $x_1$, $x_2$ cycle around a circle, we do not get
the expected (generalized) helices. This is because of our choice of the
orthonormal basis $X_1$, $X_2$ of the distribution. Changing the metric by
the choice of the orthonormal frame
\begin{align*}
X^1 &= \frac\partial{\partial x^1} - x^2\frac\partial{\partial z}
\\
X^2 &=\frac\partial{\partial x^2} + x^1\frac\partial{\partial z}
\end{align*}
the solutions $x_1(t)$ and $x_2(t)$ do not change, while $z(t)$ gets more
symmetric in the parameters:
$$
z(t) = C_4 + \frac{C_1{}^2C_2}K + \frac{C_1{}^2}Kt -
\frac1K\bigl(C_1C_3\cos(Kt+C_2)+C_1C_5\sin(Kt+C_2)\bigr)
$$
and choosing $C_3=C_4=C_5=0$ provides exactly the helices, as expected.

\begin{example}[generalized Heisenberg in 5D]
The simplest case of the previous example was at the same time the lowest
dimensional contact sub-Riemannian case. Let us look at the general
contact sub-Riemannian geometries. If $(M,g,\cD)$ is a contact sub-Riemannian manifold  of
dimensions $2n+1$, then
we can always find a local frame $X_1\dots,X_n, X_{n+1},\dots,X_{2n}$ inducing a
splitting of $\cD$ as a sum of Lagrangian subspaces (spanned by the first $n$
and second $n$ vectors) and providing the Levi
form $L$ in the canonical form, i.e. there are real positive 
functions $\la_i^2$ on $M$ 
with $L(X_i,X_{n+i})=\la_i^2$, $1\le i\le n$, $L(X_i,X_j)=0$
otherwise. Moreover, we can normalize $\la_1=1$ by the choice of the contact
form.

Let us work out one example with nontrivial function $\la=\la_2$ 
in dimension 5.
We shall deal with standard coordinates $(x^1,x^2,x^3,x^4,z)$ on 
$\Bbb R^5$ and consider $\cD$ generated by the orthonormal frame 
\begin{align*}
X_1 &= \frac\partial{\partial x^1} - 
x^3\frac\partial{\partial z}
\qquad
X_2 = \la(x^1,x^2,x^3,x^4,z)\biggl(\frac\partial{\partial x^2} - 
x^4\frac\partial{\partial z}\biggr)
\\
X_3 &= \frac\partial{\partial x^3} + 
x^1\frac\partial{\partial z}
\qquad
X_4 = \la(x^1,x^2,x^3,x^4,z)\biggl(\frac\partial{\partial x^4} + 
x^2\frac\partial{\partial z}\biggr)
\end{align*}
Proceeding exactly as in the previous example\footnote{Actually, the
computations were done with the help of the Ian Anderson's Maple package
Differential Geometry, see \begin{url}{
https://digitalcommons.usu.edu/dg/}
\end{url}, our Maple worksheet is
displayed together with the article at arxiv.}, we arrive at the following
ten equations.
Notice the terms quadratic in the functions $\al^i$, coming from the covariant derivative
(we assume $\la=\la(x^1,x^2,x^3,x^4,z)\ne0$ in all points and write $\la_z$ 
for $\frac{\partial\la}{\partial z}$, etc.).
\begin{align*}
\dot x^1 &= \al^1,\qquad \dot x^2 = \frac1\la\al^2,
\qquad
\dot x^3 = \al^3,\qquad \dot x^4 = \frac1\la\al^4,
\\
\dot z &= x^1\al^3 - x^3\al^1 + \frac1\la x^2\al^4 - \frac1\la x^4\al^2,
\\
\dot \al^1 &= - \frac{\la_{x^1}-x^3\la_z}\la(\al^2\al^2 + \al^4\al^4) - \nu\al^3,
\\
\dot \al^2 &=  \frac{\la_{x^1}-x^3\la_z}\la\al^1\al^2 +
\frac{\la_{x^3}+x_1\la_z}\la\al^2\al^3 + (\la_{x^3}+x^2\la_z)\al^2\al^4
\\&\qquad
+ (x^4\la_z - \la_{x^2})\al^4\al^4 -\la^2\nu\al^4,
\\
\dot \al^3 &= - \frac{\la_{x^3}-x^1\la_z}\la(\al^2\al^2 + \al^4\al^4) -
\nu\al^1,
\\
\dot \al^4 &=  \frac{\la_{x^1}-x^3\la_z}\la\al^1\al^4 +
\frac{\la_{x^3}+x_1\la_z}\la\al^3\al^4 - (\la_{x^4}+x^2\la_z)\al^2\al^2
\\&\qquad
- (x^4\la_z - \la_{x^2})\al^2\al^4 +\la^2\nu\al^2,
\\
\dot\nu &= \frac{2\la_z}\la(\al^2\al^2+\al^4\al^4)
.\end{align*}
In particular, if $\la_z=0$ then again $\nu$ is a free constant parameter.
These equations are again easily solved if $\la$ is a constant function.
\end{example}

\section{Canonical complements for maximally non-integrable two step geometries} \label{canny}

\newcommand{\gr}{\operatorname{gr}} 
Here we restrict to sub-Riemannian
geometries where the Levi-form (\ref{Levi})
$$
L: \cD\times \cD \to \ccD
$$
is  surjective. 
In this
case we have $TM=[\cD,\cD]$, or at the level of the associated graded
$ \gr TM=\cD\oplus [\cD,\cD]$.
Recall that with abstract indices the Levi form is denoted $L^a_{jk}$.

\newcommand{\rank}{\operatorname{rank}}

Consider the map
\begin{equation}\label{injj}
\ccD^*\otimes \cD\ni f_a^j \mapsto 
\frac{1}{2}(f_a^\ell L_{bi\ell}+f_b^\ell L_{ai\ell})L^a_{jk}h^{ik} \in \ccD^*\otimes \cD^*
.
\end{equation}
This is clearly not injective if the Levi-form $L$ is
degenerate. 
Otherwise it seems that in some broad circumstances this
map is injective, and hence is an isomorphism.  In this case we have
the following result.
\begin{thm}\label{key1}
Let $h:T^*M\to \cD$ be a sub-Riemannian geometry with surjective Levi-form $ L: \cD\times \cD \to \ccD $. If the map
(\ref{injj}) is injective then there exists a canonical metric $g$ on
$TM$ that extends $h$ on $\cD$. In particular there is such a canonical metric for contact distributions, and  also  for free distributions with $\rank (\cD)\neq 3$.
\end{thm}
\begin{proof} We have $\gr TM =\cD \oplus \ccD$, and a metric $h$ on $\cD$. 
Note that $h$ determines a metric on $\cD\wedge \cD$.  The Levi-form
determines a linear map $\cD\wedge \cD\to \ccD$ with kernel $\cN$, and
we may identify $\ccD$ with the orthogonal complement of
$\cN\subset\cD\wedge \cD $. Thus the given $h$ metric on $\cD$
determines a metric $h$ on $\ccD$.

Now we choose a splitting of (\ref{f2}) so that we have
$TM=\cD\oplus \ccD$. Clearly this splitting with the metric $h$ on
$\cD$ and $\ccD$ determines a metric on $TM=\cD\oplus \ccD$ that we
shall denote $g$. In this setting it is reasonable to write
$\cD^\perp=\ccD$ as we have identified $\ccD$ with the orthogonal
complement to $\cD$ in $TM$, with respect to $g$. Thus by our choice
of splitting we now have a metric on $TM$ and so we have the initial
data required for Theorem \ref{thm_connection}. Next we use an
adaption of part of the proof of that Theorem \ref{thm_connection}.

Recall that the symmetric 
torsion component
$T^{\cD^\perp\cD\cD^\perp}$ is determined by $g$ as in
(\ref{dpd_torsion}). In our current setting $g$ is determined by the
choice of splitting, and so the 
torsion component $T^{\cD^\perp\cD\cD^\perp}$
is entirely determined by the sub-Riemannian structure $h:T^*M\to \cD$
and the choice of splitting.

We now claim that if (\ref{injj}) is injective then  we can fix the splitting by suitably minimizing the
symmetric torsion component $T^{\cD^\perp\cD\cD^\perp}$. This result follows from Lemma \ref{Sfix} that follows.

The final statement of the Theorem now follows from Lemma
\ref{fifteen} and Lemma \ref{sixteen}.
\end{proof}

For Lemma \ref{Sfix}  we need a preliminary result:
\begin{lemma}\label{tform}
We consider the torsion component $T^{\cD^\perp\cD\cD^\perp}$ in the setting of the Theorem \ref{key1}. Under a change of splitting of the sequence (\ref{f2}) given by $f_a^j\in \ccD^*\otimes \cD$ this transforms according to 
$$
T_{bai}\mapsto T_{bai} + \frac{1}{2}(f_a^jL_{bij}+f_b^jL_{aij}).
$$ 
\end{lemma}
\begin{proof} In order to facilitate the comparison we view the torsion component $T^{\cD^\perp\cD\cD^\perp}$ and its analogue in the new splitting
each as a tensor in $\ccD^*\otimes \cD^*\otimes \ccD$. 

We start with some initial splitting that defines $g$. In this splitting
$\cD^\perp$ is the image of $\ccD$, and we may identify these two and
view the change of splitting as a map $f:\cD^\perp\to \cD$.  We adorn
with hats the objects in the new splitting. For example $\widehat{g}$
is the metric determined by the new splitting. Then for
$\tilde{X},\tilde{Y}\in \ccD\subset \gr(TM) $ with representatives $
X,Y \in \cD^\perp\subset TM$, and $\widehat{X}, \widehat{Y}\in
\widehat{\cD}^\perp\subset TM$ respectively, we have
$$
\widehat{X} = X+f(X), \qquad \mbox{and}\qquad \widehat{Y} = Y+f(Y)
$$
in $TM$, 
and so 
$$ \langle \tilde{X},\tilde{Y} \rangle_h = \langle X,Y \rangle_g =
\langle \widehat{X},\widehat{Y} \rangle_{\widehat{g}}= \langle X
+f(X),Y+f(Y) \rangle_{\widehat{g}},
$$ since the metrics $g$ and $\widehat{g}$ are each compatible with
the metric $h$ on $\gr (TM)$ via the respective splittings.

Now let $Z\in \cD$. Note that  $\cD\subset TM $ is fixed in the change of splitting so we have $\widehat{Z}=Z$.  
 We have, using the 
formula (\ref{diag_dpd}), 
$$
\begin{aligned}
\langle \widehat{X}, \widehat{T}(\widehat{Z}, \widehat{X})\rangle_{\widehat{g}} & = 
\frac{1}{2} Z\cdot ||X+ f(X) ||_{\widehat{g}}^2 + \langle X+ f(X),  [X+ f(X), Z]_{\widehat{\cD}^\perp} \rangle_{\widehat{g}}.
 \end{aligned}
$$
 Thus, using the observation in the previous display, and again formula
(\ref{diag_dpd}), this simplifies to
$$
\begin{aligned}
\langle \widehat{X}, \widehat{T}(\widehat{Z}, \widehat{X})\rangle_{\widehat{g}}
 & = 
\frac{1}{2} Z\cdot ||X ||_g^2 + \langle X,  [X+ f(X), Z]_{\cD^\perp} \rangle_g\\
& = \langle X, T (Z,X)  \rangle_g + \langle X, [f(X),Z]_{\cD^\perp} \rangle_g .
 \end{aligned}
$$
That is (after multiplying through with $-1$)
$$ 
\tilde{X}_b  \widehat{T}^b{}_{ai}\tilde{Z}^i\tilde{X}^a  = \tilde{X}_b  T^b{}_{ai}\tilde{Z}^i\tilde{X}^a + \tilde{X}_b L^b_{ij} f^j_a\tilde{X}^a\tilde{Z}^i .
$$
\end{proof}

\begin{lemma}\label{Sfix}
In the setting of the Theorem, and if (\ref{injj}) injective, there
exists a unique splitting of (\ref{f2}) such that the torsion component
$T^{\cD^\perp\cD\cD^\perp}$ satisfies
$$
T^b{}_{ai}L^a_{jk}h^{ik}=0,
$$
where $h$ is the canonical metric $h$ on $\ccD$ (as discussed above).
\end{lemma}
\begin{proof}
Fix a choice of splitting and denote by $\tilde{g}$ the resulting
metric on $TM$. With respect to $\tilde{g}$ we have the torsion component $T^b{}_{ai}=T^{\cD^\perp\cD\cD^\perp}$, and we recall that this is
  symmetric. We consider 
$U^b_j:=T^b{}_{ai}L^a_{jk}h^{ik}\in \ccD\otimes \cD^*$.

A change of splitting of (\ref{f2}) is a map $f:\ccD\stackrel{\tilde{g}}{=}\tilde{\cD}^\perp\to \cD$,
that in abstract indices we denote $f_a^i$. As established in
Lemma \ref{tform}, under such a change we have 
$$
T_{bai}\mapsto T_{bai} + \frac{1}{2}(f_a^jL_{bij}+f_b^jL_{aij}),
$$ 
where indices were lowered using the canonical metric $h$ on $\ccD$.
Thus the result follows immediately from the fact 
that the map (\ref{injj})
is an isomorphism.
\end{proof}

\begin{lemma}
  The map (\ref{injj}) is injective if and only if the map
\begin{equation}\label{need}
 (\ccD)^*\otimes \cD\ni f_a^j \mapsto 
\frac{1}{2}(f_a^\ell L_{bi\ell}+f_b^\ell L_{ai\ell}) \in S^2\ccD^*\otimes \cD^*
  \end{equation}
is injective.
\end{lemma}
\begin{proof}{} 
If $f_a^\ell$ is in the kernel of (\ref{need}) then it in the kernel
of (\ref{injj}), as the latter is a composition of the map (\ref{need})
with a subsequent map.
  
For the forward implication we suppose that $f_a^\ell$ is in the kernel of
(\ref{injj}). Then the right hand side of (\ref{injj}) is zero and by
contracting in $f^{jb}$ we see
  $$
0=||U ||^2=U_{abi}U_{cdk}h^{ik}g^{ac}g^{bd} ,
  $$
where
$$
U_{abi}:=\frac12( f^\ell_a L_{bi\ell}+ f^\ell_a L_{bi\ell}) ,
$$
and indices have been lowered (and raised) using the metric
$h_{ab}$ on $\mathcal{Q}$  (and its inverse). Thus $U=0$ and $f_a^\ell$ is in the kernel of (\ref{need}). 
  \end{proof}
\begin{lemma}\label{fifteen}
For contact sub-Riemannian geometries (\ref{need}) is injective.
\end{lemma}
\begin{proof}
  In this case the map is equivalent to
  $$
 \cD\ni f^j \mapsto 
f^\ell L_{i\ell} \in  \cD^*
$$
which is injective since the Levi form $L_{ij}$ is non-degenerate.
  \end{proof}
\begin{lemma}\label{sixteen}
For free sub-Riemannian geometries (\ref{need}) is injective if and only if $n\neq 3$.
\end{lemma}
\begin{proof}
First observe that injectivity for the case $n=\dim (\cD)=2$ is covered by Lemma
\ref{fifteen}. Thus we now consider only $\dim (\cD)\geq 3$. 

In any free case we may first identify $\ccD$ with $\Lambda^2\cD$ by
  using $L^a_{ij}$.  With this done, the map (\ref{need}) becomes
  \begin{equation}\label{helps}
  \cD\otimes \Lambda^2\cD\ni f^{jmn}
  \mapsto
  \frac{1}{4}(\delta^p_i f^{qmn}-\delta^q_i f^{pmn}+ \delta^m_i f^{npq}-\delta^n_i f^{mpq})
  \in S^2(\Lambda^2 \ccD) \otimes \cD^* ,
  \end{equation}
  where we have used $h$ to raise and lower indices (and identify
  $\cD$ with $\cD^*$).

 Now we consider any case  $\dim (\cD)\geq 4$ and  assume  that  $f^{jmn}$ is in the kernel of the map (\ref{helps}).
  It follows that
  \begin{equation}\label{key4}
\delta^{[p}_i f^{qmn]}=0 ,
  \end{equation}
  where (here and below) the notation $[\cdots]$ indicates the
  completely skew part indicated by the  the enclosed indices.
Then from (\ref{key4}) it follows that 
$$
f^{[qmn]}=0.
$$
Next contracting the right hand side of (\ref{helps}) with $\delta^i_p$ gives
$$
(n-1)f^{qmn}+f^{nmq}-f^{mnq}=0 ,
$$
which, with the previous display, gives
$$
f^{qmn}=0. 
$$
Thus (\ref{helps}), equivalently (\ref{need}), is injective if
$n\geq 4$. 

Now consider the case  that $\dim (\cD)= 3$. Suppose that
$f^{jmn}\in \Lambda^3\cD$. Then the image of (\ref{helps})
$$ \frac{1}{4}(\delta^p_i f^{qmn}-\delta^q_i f^{pmn}+ \delta^m_i
f^{npq}-\delta^n_i f^{mpq})
$$ is exactly $ \delta^{[p}_i f^{qmn]}$ which vanishes, as $\dim
(\cD)<4$. Thus $ \Lambda^3\cD$ is in the kernel of (\ref{helps}) if
$\dim (\cD)= 3$. (In fact it is straightforward to show that $
\Lambda^3\cD$ is exactly the kernel of (\ref{helps}), but we do not
need that here.) This completes the proof.
\end{proof}

\end{document}